\documentclass[11 pt]{article}
\usepackage[margin=1.5in]{geometry}
\usepackage{amsmath}
\usepackage{amssymb}
\usepackage{amsthm}
\usepackage{amsfonts}
\usepackage[none]{hyphenat}
\usepackage{fancyhdr}
\usepackage{graphicx}
\usepackage{float}
\usepackage[nottoc,notlot, notlof]{tocbibind}

\usepackage[varg]{txfonts}
\usepackage{color}

\definecolor{blue}{rgb}{0,0,1}
\definecolor{red}{rgb}{1,0,0}
\definecolor{green}{rgb}{0,.6,.2}
\definecolor{purple}{rgb}{1,0,1}
\definecolor{brown}{rgb}{.59,.29,0}

\long\def\red#1\endred{{\color{red}#1}}
\long\def\blue#1\endblue{{\color{blue}#1}}
\long\def\purple#1\endpurple{{\color{purple} #1}}
\long\def\green#1\endgreen{{\color{green}#1}}
\long\def\brown#1\endbrown{{\color{brown}#1}}

\begin{document}

\begin{titlepage}
\end{titlepage}

\title{\bf{Integral Solutions of :$x^p-my^p=zw.$}}
\author{\bf{  ANTENOR OSTINE  }}

\maketitle

\begin{abstract}

In this paper, using a basic algebraic approach, we demonstrate  that, 
if the integers $x,y$ and $z$ are pairwise relatively prime, and 
$x^p-my^p=zw$ where $m$ is an integer and $p$ is a prime, then  we can write $x,y,m,z$ and $w$ as explicit 
expressions in 7 integers.

\end{abstract}
\newpage
\tableofcontents
\thispagestyle{empty}
\clearpage
\setcounter{page}{1}

\section{Introduction}

{\bf{ }}           
Diophantus of Alexandria was the first to prove that, in algebra, the set of all sum of two squares is closed under multiplication:
$$  (a^2+q^2)( b^2+r^2)=(ab \pm qr  )^2+(ar\mp bq)^2$$ 
Brahmagupta (598-668) later proved a more general identity which showed that\cite{DBHS2}:
$$   (a^2-mq^2)( b^2-mr^2)=(ab \mp mqr  )^2-m(ar\pm bq)^2    $$
In 1920, LE Dickson built upon this identity and proved that all integral solutions of the equation: $$ x^2-my^2=zw$$ are given by:
$$\pm x= eln+fnq-flr-gqr$$
$$y=nq+lr$$
$$m=f^2-eg$$
$$z=el^2+2flq+gq^2$$
$$w=en^2-2fnr+gr^2$$ \cite{DBHS1}\\
In our approach, we expand even further LE Dickson's findings, and describe the integral solutions of the equation $x^p-my^p=zw$ as expressions of the integers $e,f,g,l,q,n,$ and $r$.

\section{ Methodology} 
First, we define a linear equation: $qx-uy+rz=0$ where the integral triplet $(x,y,z)$ is pairwise coprime, and $\gcd(q,u)=1$. Then, we determine its intersection with the equation of the form $x^p-my^p=zw$ where $m,w$ are integers and $p$ is a prime. 

\section{ Integral solutions of $x^p-my^p=zw$.}    

{\bf{Theorem 1.1}.{\it{ We fix a prime $p$. Let $x,y,z,m,w$ be non-zero integers such that $x$, $y$ and $z$ are pairwise coprime and such that $x^p-m y^p=zw$. Then there are integers $e,f,g,l,q,n,r$ such that : 
\begin{eqnarray}
x & = & eln-(\sum_{k=1}^{p-1}{ p\choose k}e^{p-k-1}l^{p-k}f^kq^{k-1}+gq^{p-1})r+fy\\
y &=& nq+e^{p-2}l^{p-1}r\\
m &= &f^p-eg\\
z& =& \sum_{k=0}^{p-1}{ p\choose k}e^{p-k-1}l^{p-k}(fq)^k+gq^p\\
w& = & \frac{1}{q^p}\Bigl [\sum_{k=0}^{p-1}{ p\choose k}z^{p-k-1}(-r)^{p-k}((el+fq)y)^k+ey^p\Bigr]
\end{eqnarray} 
where $q\neq 0$, and $(e,q)=(l,q)=(n, r)=1.$ }}}\\

{\it{Proof}}.\\
Let $x,y,z$ be 3 pairwise coprime nonzero integers and  $p$ be a prime such that:
\begin{eqnarray}
x^p-my^p=zw
\end{eqnarray}
Bezout's identity implies, there exist infinitely many integers $a,b,c,d$ such that:
\begin{eqnarray}
ax-bz=cy-dz=1
\end{eqnarray}

Not all of those integral coefficients are equal to zero. For instance, if $a$ and $b$ are be both equal to zero, then $0=1$, which is absurd. For similar reason, we cannot have either  $c=d=0$. Thus, at least, one of them must be non zero. \\Note that if $z\neq 0$, then any solution $(a,b)$ can be replaced by $(a+z,b+x)$. This way, $a$ can be arranged to be non-zero. If $z=0$, then $x=\pm 1$ by coprimality, and we can take $a=x\neq0$.  For $c$ and $d$, one may argue similarly.  So we can take $a, b, c, d $ satisfying (7), with $a\neq 0$ and $c\neq 0$.\\ If a prime divides $a$ and $b$, then it must divide 1, which is absurd since no prime divides 1. The same argument holds if a prime divides both $c$ and $d$. Then, we have $\gcd(a,b)=\gcd(c,d)=1$. \\

From equation (7), we obtain:
\begin{eqnarray}
ax=-z(d-b)+cy
\end{eqnarray}
Put  $\gcd(a,c)=h$. There exist coprime integers $u$ and $q$ so that $a=qh$ and $c=uh$. Let a prime divide $h$ and $z$. Then, according to (7) it must also divide 1, which is a contradiction since no prime divides 1. Thus, $h$ divides $d-b$. There is an integer $r$ such that $d-b=rh$. It follows that
\begin{eqnarray}
qhx=-zhr+uhy
\end{eqnarray}

Since $a=qh\neq 0$, then $h\neq 0$. We have:
\begin{eqnarray}
qx=-zr+uy
\end{eqnarray}
 
 We multiply both sides of equation (6) by $q^p$. Then we substitute $qx$, and obtain:
\begin{eqnarray}
(-zr+uy)^p-m(qy)^p=zwq^p
\end{eqnarray}
Using the binomial theorem, we expand the lefthand side,
\begin{eqnarray}
\sum_{k=0}^{p}{ p\choose k}(-zr)^{p-k}(uy)^k-m(qy)^p=zwq^p
\end{eqnarray}
Which gives us after simple algebraic manipulations:
\begin{eqnarray}
\sum_{k=0}^{p-1}{ p\choose k}(-zr)^{p-k}(uy)^k+(uy)^p-m(qy)^p=zwq^p\\
\sum_{k=0}^{p-1}{ p\choose k}(-zr)^{p-k}(uy)^k+y^p(u^p-mq^p)=zwq^p
\end{eqnarray}
The righthand side of (14) is divisible by $z$, and hence so is the lefthand side. All terms in the sum over $k$ are divisible by $z$. Hence, the second term in the lefthand side is also divisible by $z$. Since $\gcd( y,z)=1$, thus there exists an integer $e$ such that:
\begin{eqnarray}
u^p-mq^p=ze\\
\sum_{k=0}^{p-1}{ p\choose k}(z)^{p-k-1}(-r)^{p-k}(uy)^k+y^pe=wq^p
\end{eqnarray}
We have shown between (8)-(9) that $u$ and $q$ are coprime. Let a prime divide both $e$ and $q$, thus according to (15), it must divide $u$ in contradiction to the fact that $\gcd(u,q)=1$. Hence, $e$ and $q$ are coprime. Thus, it's possible to find 2 integers $l,f$ such that:
\begin{eqnarray}
u=el+fq
\end{eqnarray}

By substituting $u$ in equation (15), we obtain:
\begin{eqnarray}
(el+fq)^p-mq^p=ze\\
\sum_{k=0}^{p}{ p\choose k}(el)^{p-k}(fq)^k-mq^p=ze\\
\sum_{k=0}^{p-1}{ p\choose k}(el)^{p-k}(fq)^k+f^pq^p-mq^p=ze\\
\sum_{k=0}^{p-1}{ p\choose k}e^{p-k}l^{p-k}(fq)^k+(f^p-m)q^p=ze
\end{eqnarray}
The righthand side of (21) is a multiple of $e$, and hence so is the lefthand side. All terms in the sum over $k$ are a multiple of $e$. Hence, the second term in the lefthand side is also a multiple of  $e$. Since $\gcd( e,q)=1$, thus there exists an integer $g$ such that:
\begin{eqnarray}
f^p-m=eg
\end{eqnarray}
It follows that:
\begin{eqnarray}
\sum_{k=0}^{p-1}{ p\choose k}e^{p-k-1}l^{p-k}(fq)^k+gq^p=z
\end{eqnarray}
Substituting $z$ and $u=el+fq$ in equation (10) yields:
\begin{eqnarray}
xq=-(\sum_{k=0}^{p-1}{ p\choose k}e^{p-k-1}l^{p-k}(fq)^k+gq^p)r+y(el+fq)\\
xq=-e^{p-1}l^pr-(\sum_{k=1}^{p-1}{ p\choose k}e^{p-k-1}l^{p-k}(fq)^k+gq^p)r+y(el+fq)\\
xq=el(y-e^{p-2}l^{p-1}r)-(\sum_{k=1}^{p-1}{ p\choose k}e^{p-k-1}l^{p-k}(fq)^k+gq^p)r+fyq
\end{eqnarray}

According to (17), if a prime divides both $q$ and $el$, then it must divide $u$. This is absurd since $u,q$ are relatively prime. Hence, $(e,q)=(l,q)=1$. \\
The lefthand side of (26) is divisible by $q\neq 0$, and hence so is the righthand side. All terms in the sum over $k$ and $gq^p$ are divisible by $q$. Hence, the first term in the righthand side is also divisible by $q$. Since $q$ and $el$ are relatively prime, then $q$ must divide $y-e^{p-2}l^{p-1}r$. It follows there exists an integer $n$ such that:
\begin{eqnarray}
y=nq+e^{p-2}l^{p-1}r
\end{eqnarray}
Which gives us: 
\begin{eqnarray}
x=eln-(\sum_{k=1}^{p-1}{ p\choose k}e^{p-k-1}l^{p-k}f^kq^{k-1}+gq^{p-1})r+fy
\end{eqnarray}
Equation (16) becomes after subsituting $u=el+fq$:
\begin{eqnarray}
q^pw=\sum_{k=0}^{p-1}{ p\choose k}(z)^{p-k-1}(-r)^{p-k}(y(el+fq))^k+ey^p
\end{eqnarray}
where $q\neq 0$ since $a=hq\neq 0$.\\
It is worth noting that both sides of all the equations from (11) through (16) but (15) are divisible by $q^p$.  Therefore, after dividing both sides of (29) by $q^p$, we obtain:
\begin{eqnarray}
w= \frac{1}{q^p}\Bigl [\sum_{k=0}^{p-1}{ p\choose k}z^{p-k-1}(-r)^{p-k}((el+fq)y)^k+ey^p\Bigr]
\end{eqnarray}
where $w$ is an integer and $q\neq 0$.\\
 At this point, (1)-(5) have been shown in equations (22), (23), (27) (28), and (30). The coprimality of $(e,q)$ and $ (l,q)$ were respectively shown between equations(16)-(17) and  (26)-(27).\\ Let a prime divide both $n$ and $r$. Then it also divides $x$ and $y$ according to equations (27)-(28). This is a contradiction since $\gcd(x,y)=1$. Hence, $n$ and $r$ are coprime. \\$\square $\\

{\bf{Corrolary 1.1}{\it{ We fix a prime $p$. Let $x,y,z,m,w$ be nonzero integers such that $x$, $y$ and $z$ are pairwise coprime and such that $x^p-m y^p=zw$. Then one can obtain explicit expressions for $(x,y,z,m,w)$ in terms of $(e,f,g,l,q,n,r).$ }}}\\
 
{\it{Proof}}.\\
 Substitute the expression for $y$ in (2) into (1) and (5), and the expression for $z$ in (4) into (5).$\square$
 
\section{Conclusions}
We demonstrated, when they exist, it is possible to express the integral solutions of the diophantine equation of the form $x^p-my^p=zw$ as explicit functions of 7 integers. This result can be helpful in identifying whether a specific integral solution set of this equation exists. For if it does, then we should be able to find a corresponding set of integers $e,f,g,l,q,n,r$ for which that assumption holds. 

 \section{Acknowledgements}
I would like to thank Dr. Roelof Bruggeman for his unwavering patience and invaluable feedbacks.

\pagebreak


\begin{thebibliography}{ }




\bibitem{DBHS1}
Dickson, L.E., \textit{Integral solutions of $x^2-my^2=zw$.}  Bull. Amer. Math. Soc. 29 (1923), no. 10, 464--467. Retrieved from https://projecteuclid.org/euclid.bams/1183485743 .

\bibitem{DBHS2}Wikipedia contributors. (2019, May 8). Brahmagupta–Fibonacci identity. In Wikipedia, The Free Encyclopedia. Retrieved 20:09, March 16, 2020, from https://en.wikipedia.org/w/index.php?title=Brahmagupta.


\end{thebibliography}
\end{document}